\documentstyle{article}
\input epsf

\newcommand {\Si}{\Sigma}

\newcommand{\R }{I \!\! R}
\newcommand{\Z }{Z \!\!\! Z}

\newcommand{\D}{ {\cal D}   }
\newcommand{\Dp}{{\cal D ^{\perp}}}

\newtheorem{proposition}{Proposition}
\newtheorem{theorem}{Theorem}
\newtheorem{lemma}{Lemma}

\title{\sc  Characteristic classes for the degenerations of  2-plane
fields in four dimensions}

\author{ Maxim Kazarian \\ Steklov Mathematical Institute,  42 Vavilova
St., 117966,\\ Moscow GSP-1, Russia, e-mail: kazarian@ium.ips.ras.ru \\
Richard Montgomery \\ Mathematics Dept. UCSC, Santa
Cruz, CA 95064 \\USA, e-mail:  rmont@cats.ucsc.edu\\
Boris Shapiro \\Dept of Mathematics, Univ. of Stockholm, S-10691,\\
Stockholm, Sweden, e-mail: shapiro@matematik.su.se\\
\\}
\begin{document}

\maketitle

{\sc abstract} { \bf Given a distribution of $k$-planes on a
 manifold,   consider the degeneration locus $\Sigma_I$ consisting of
points  where  the distribution
has Lie-bracket growth  vector less than or equal $I$, a fixed integer
vector. We calculate the characteristic classes  associated to the
$\Sigma_I$ for generic $2$-plane distribution  on a 4-manifold.}

\section{Results and Background}

\subsection{Generalities, Setting and Results.}

A distribution $\D$ of k-planes on an n-dimensional manifold
$Q$ can be thought of as either a subbundle
$\D \subset TQ$ of the tangent bundle or as a locally free
sheaf of smooth vector fields.  We use the same notation for
both. Write $\D^2 = \D + [\D, \D]$ and more generally
$\D^{j+1} = \D^j + [\D, \D^j]$.  These are sheaves
of modules  of
vector fields (over the ring of smooth functions).
We are interested in distributions
such that for $r$ large enough we obtain
all vector fields by this procedure:
$$\D^r = T$$
where $T$ denotes the sheaf of all vector fields.
These are called
{\sc completely nonholonomic}
distributions.
We thus have a filtration
$$\D \subset \D^2 \subset \ldots \subset \D^r =  T$$
by subsheaves of the sheaf of all vector fields.
Write $\D^j (q) \subset T_q Q$ for the vector subspace
obtained by evaluating all vector fields in $\D^j$ at the point
$q \in Q$, and
$$n_j (q) = dim(\D^j (q)).$$
The first $r$
such that $\D^r (q) = T_q Q$
is called the degree of nonholonomy at
$q$ and the
nondecreasing list of dimensions
$$I(q)=(n_1(q) , n_2 (q), \ldots,
n_r(q))$$
 is called the {\sc growth vector} at $q$.
If the $n_j (q)$ are constant in a neighborhood
of  $q$ then the $\D^j$
correspond to vector bundles of this rank
(in this neighborhood) which we will also denote by $\D^j$.
In this case,
we say that
$q$ is a {\sc regular point} for the distribution.
Otherwise, some of the
ranks $n_j$ jump as we pass through $q$
and we have to think of the corresponding $\D^j$ as
sheaves.

{\sc Example.}  A contact distribution has
growth vector $(k, k +1)$ with $k$ even.

{\sc Definition.} An {\sc Engel distribution} is
a rank two  distribution on a 4-dimensional manifold
whose growth vector is $(2,3,4)$ everywhere.

What makes Engel distributions remarkable is that they are
topologically stable, and topologically 
stable distributions are quite rare,   
occuring only in dimensions $(k,n) = (1,n), (n-1,n)$ and $(2,4)$. 
The only stable
regular distributions are the line fields, the contact fields,
an odd-rank analogue of contact (sometimes called
pseudo-contact)  and the Engel
distribution. 
( See the next section, for the definition of
stability, and some details.)

If one slightly perturbs any given
distribution  of 2-plane fields on a  4-manifold then it will become
Engel on an  open dense subset (\cite {Zhit},  \cite{Gersh}).  On the
other hand (see propositon 1 below) if an oriented 4-manifold admits an
oriented Engel distribution, then that manifold is
 parallelizable. {\bf  Thus there are topological obstructions
to making the 2-plane field globally  Engel.  Our goal is to 
understand these obstructions.}

A basic notion in this endeavour will be the degeneraion
locus of a distribution. We will first need to 
establish a partial order for growth vectors.  
Declare that $J = (m_1, m_2,  , m_s = n)  \le I = (n_1, n_2, \ldots
, n_r = n)$ if and only if  $m_i \le n_i$ for  $i = 1, \ldots , r$.
Note that the two vectors may have  a different number of components.
For fixed $k = n_1$ and $n$ there is exactly one maximal growth vectors.
Its components,
except for possibly
the last one $n_r = n$ are the dimensions of the subspaces of the standard grading
of the free Lie algebra on k elements.  For typical distributions
the growth vector will be maximal at most points of the manifold.

{\sc Definition.} The degeneration locus of a distribution
$\D$ on a manifold $Q$ is the set of all points $\Sigma = \Sigma(\D) \subset
Q $ whose growth vector is less than maximal.  The degeneration locus of
type $I$,   denoted $\Sigma_I = \Sigma_I(\D)\subset Q$, is  the subset of
all points $q\in Q$ at which   the growth  vector is less than or equal $I$.

The Thom transversality theorem  implies that for typical $\D$ 
 all of the $\Sigma_I$ are nice subvarieties. They
stratify  the manifold.  
 
The Engel growth vector $(2,3,4)$ 
is the maximal growth vector
for a rank two distribution in 4-space.  
The smaller 
growth vectors $I_1 = (2,2,4)$, $I_2 = (2,3,3,4)$ are not ordered relative
to each other.  All other growth vectors
are dominated by these two.  Consequently:
 $$\Sigma = \Sigma_1 \cup \Sigma_2$$
The following alternative
descriptions of the $\Sigma_i$ may be more transparent: 
$$\Sigma_1  = \{ q : dim(\D^2 (q)) \leq 2 \}$$
$$\Sigma_2= \{ q : dim(\D^3 (q)) \leq  3 \}.$$
Zhitomirskii has shown that generically each piece $\Sigma_i$
is a smooth 2-dimensional surface.  Our question becomes, 
what is the topological meaning of the degeneration locus $\Sigma$
and its pieces $\Sigma_1$ and $\Sigma_2$?  
The situation is complicated by the fact that
when the two pieces do intersect, they never do
so transversally, but rather along a curve, denoted
by C below.

In general, the condition that the
growth vector $I(q)$ is less than or equal to  some $I=(n_1,...,n_r)$
defines  a natural $Diff_Q$-invariant subset $\cal M_I$ in the space of
$r$-jets of  sections of the bundle $G_{k,n}(TQ)\to Q$ whose
fibers
$G_{k,n}(T_q Q)$ are the  Grassmannians of k-planes in the
n-dimensional tangent spaces  $T_qQ$, 
and where $Diff_Q$ denotes the group of diffeomorphisms of $Q$.
$\Sigma_I(\D)$ is  the pullback of the
intersection of the $r$-jet  extension of $\D$ with $\cal
M_I$.  When we say `typical' and `generic' for $\D$, this means that the
$r$-jet extension of $\D$ is transversal to the $\cal M_I$.   By results of
Thom, (see \cite {HaefligerKos}) for a generic $\D$ there exists a
universal formula for the Stiefel-Whitney and, in the  appropriate setting,
the Chern classes
 which are Poincare dual to $\Sigma_I(\D)$ in terms of the characteristic
classes of $\D$ and $Q$.  In this paper we obtain this
formula for the Engel degenerations. 
The main result of this note is as follows.

\begin{theorem}
Let $\D$ be a real oriented 2-plane field on a closed
oriented 4-manifold $Q$. Then the obstructions to $\D$  being
 Engel   are the  Chern classes $c_1(\D^\perp)$
and $c_1(\D^*)$. Here  $\D^\perp=TQ/\D$ and $\D^*$ denotes the dual bundle.
For generic $\D$ the
 first class is represented by $\Sigma_1 \subset \Sigma$  which
is a smooth 2-surface (except for possibly a finite
number of points) 
with a global co-orientation.  
 The second class is represented by
$\Si$ with a certain co-orientation
on $\Si \setminus C$ where $C =\Si_1\cap\Si_2$.
The 1-cycle C is homologous to
zero in $\Sigma_1$.
\end{theorem}

{\sc Remark.} $\Si_1$ is canonically
co-oriented as the class representing $c_1(\D^\perp)$.
In fact it is the zero locus of a section
of $\Dp$ which is transverse to the zero section.
The co-orientation which $\Si_1$ receives. as
part of the cocycle representing $\D^*$
reverses relative to the canonical
co-orientation as the curve C is crossed.
$\Si_2$ has no canonical
coorientation and may be nonorientable.

 If the distribution $\D$ or the underlying manifold $Q$ are
not oriented then
the degeneration loci  will not yield well-define integer 
homology classes,
but rather only classes $mod2$. Let $w_i(\xi)$ denote the ith
Steiffel-Whitney class of a real vector bundle, and write $w_i(Q) =
w_i(TQ)$.  

\begin{theorem}
The $\Z_2$-homology classes of $\Si_1$, $\Si_2$, $\Si_1\cup\Si_2$ and
 $\Si_1\cap\Si_2$ are
dual to
$w_1^2(D)+w_2(D)+w_2(Q)$, $w_1^2(D)+w_1(Q)w_1(D)+w_1^2(Q)+w_2(Q)$,
 $w_2(D)+w_1^2(Q)+w_1(D)w_1(Q)$
 and
$w_1(D)(w_1^2(D)+w_2(D)+w_2(Q))$ respectively.
\end{theorem}

Two of us (M. Kazarian and B. Shapiro),
have found 
corresponding duality formulae for the case of generic distributions of
arbitrary rank k in n dimensions, and plan to publish
this in a future paper.  
 
\subsection{Acknowledgements} 
We would like to thank Ivan Kupka for a
careful reading and critique of
an earlier version of this paper. R. Montgomery would like
to thank the NSF, grant number DMS 9400515 and
the University of California at Santa Cruz, FRC grants
for support.  M. Kazarian would like to thank
the Russian Foundation
for Basic Researches (RFSB) grant 95-01-01122a for support.  
  
\subsection{Properties of rank two distributions on
4-manifolds}

We  formulate some  known properties of rank
two distributions $\D$ on a 4-manifold $Q$.

An oriented 4-manifold $Q$ admits an oriented
2-plane  distribution $\D$ if and only if $\chi[Q]\equiv 0$ mod 2 and
$\chi[Q]=\tau  [Q]$ mod 4.  Here $\chi$ is the Euler class and
$\tau$ the signature.   This condition is equivalent to requiring
that
the manifold admit two almost complex structures, one consistent  with
the given orientation, and the other consistent with the opposite
orientation.  (See \cite{Atiyah} and  also \cite{Matsushita} and
references therein.)

If  the distribution $\D$  has
 rank two then it is locally spanned by  two non-vanishing vector
fields $X$ and $W$.  Locally:
$$\D^2 = span \{X, W, [X,W] \}$$
and
$$\D^3 = span \{ X, W, [X,W], [X,[X,W]], [W,[X,W]] \}.$$
It follows that  $\D$ is Engel if and only if for some functions
$a$, $b$ the vector fields
 $X, W, [X,W]$ and $a[X,[X,W]] + b[W,[X,W]]$
form a basis for the tangent space.

{\bf Engel's Theorem.}  (See \cite{VG}, \cite{EDS}, p. 50. )
{\it If $\D$   is Engel at the point $q \in Q$ then it admits a local
frame $X$ and $W$ and $Q$ admits coordinates
$(x,y,z, w)$ centered at  $q$ such  that
$W = {\partial \over {\partial w}}$
and
$X = {\partial \over {\partial x}} +
 w{\partial \over {\partial y}}
+ y{\partial \over {\partial z}}$.}

To clarify the special nature of Engel
distributions, we   recall the notion of ``stability''
in singularity theory. A distribution
germ is called {\sc stable}  if every sufficiently nearby
distribution germ is  diffeomorphic to it, sufficiently near
being measured in the $C^r$-topology for some $r$.

{\bf Stability Theorem}  (See \cite{VG} and also \cite{Mont}.)
{\it The only stable distribution germs
occur in dimensions $(k,n) = (1,n), (n-1,n)$ or $(2,4)$.
In each case there is a unique stable regular representative.
These are the line fields, contact (or even-contact) fields,
and the Engel distributions.}

The generic degenerations of Engel structures
were classified, up to codimension 3, by Zhitomirskii \cite{Zhit}.
The first degeneracies occur in codimension 2.  There are two of them,
one for $\Sigma_1$, and one for $\Sigma_2$.  Both 
are stable.  Corresponding normal
forms were found by Zhitomirskii
and are repeated below. The main
result of Zhitomirskii's which we will
be using is that 
{\sc for a Whitney open
and dense set of rank two distributions
on a 4-manifold, the degeneration loci $\Sigma_1 $ and $\Sigma_2 $
 are smooth 2-dimensional embedded surfaces.} 
Outside the intersection of $\Sigma_1$ and $\Sigma_2$,
this  follows from Thom's transversality theorem
(\cite{Arn}, p. 38, and references therein).

We will also require the normal forms and higher
singularities however.  There are two types
of codimension 3 degenerations.  One occurs
along $C = \Si_1 \cap \Si_2$ which is
typically a curve when nonempty. The
other occurs along a curve L in $\Si_2$
which separates the degenerating Engel line
field near $\Si_2$ into hyperbolic and
elliptic degenerations.  The corresponding
normal forms are given in terms of 1-forms
$\omega_1, \omega_2$ which frame the forms
annihilating $\D$. They are
\begin{equation}
\Si_1 \;: \omega_1 = dx_1 +x_3 ^2 dx_4 , \omega_2 = dx_3 + x_3 x_4 dx_4
\label{eq:Z1}
\end{equation}
\begin{equation}
\Si_2 \;: \omega_1 = dx_1 +x_3 dx_4 , \omega_2 = dx_3 + 
{1 \over 3}(x_3 ^2 + x_3 x_4 )dx_4
\label{eq:Z2A}
\end{equation}
or
\begin{equation}
\Si_2 \;: \omega_1 = dx_1 +x_3 dx_4 , \omega_2 = dx_3 + x_3 ^2 x_4 dx_4
\label{eq:Z2B}
\end{equation}
The two normal  forms for $\Si_2$ correspond to
the two regions of $\Si_2$ separated by the curve L.
For the points of C and L there is a single functional
modulus f.  Their normal forms are:
\begin{equation}
C \;: \omega_1 = dx_1 +x_3 ^2 dx_4 , \omega_2 = dx_3 + x_3 x_4 dx_4
\label{eq:Z3}
\end{equation}
\begin{equation}
L \;: \omega_1 = dx_1 +x_3 ^2 dx_4 , \omega_2 = dx_3 + x_3 x_4 dx_4
\
\label{eq:Z4}
\end{equation}
There are also codimension 4 singularities
which occur along the curves C and L as isolated
points.  Fortunately we will not need explicit
normal forms for them. (None have been calculated!)

We re-iterate the main picture.
$\Sigma$ is a smooth 2-dimensional surface
away from the curve C. Near C it  is diffeomorphic
to the interstion of two 2-planes in 4-space
which intersect on a line.  {\bf See figure 1 below.}

\subsection{ The Engel line field and parallelizability}

We recall  why
Engel manifolds are parallelizable, up to orientation problems.

It will help to
recall the basic Lie algebra structure associated
to a distribution.  Let $[,]$ denote the operation
of Lie bracket of vector fields.  If
 $f$ and $g$ are smooth functions on an
open set $U \subset Q$ and $X \in \D^j (U)$, $Y \in \D^k (U)$
then
$$[fX, gY] = fg [X,Y] mod(\D^{l})(U)$$
where $l = max(j,k)$.  In
other words, Lie bracket induces
well-defined maps $\D^j \otimes \D^k \to T/\D^l$
where $T$ denotes the (sheaf of the) tangent bundle.

Associate to our filtration
$\D \subset \D^2 \subset ...\subset T$
  the corresponding graded object:
$$Gr(T) := \D \oplus V_2 \oplus V_3 \oplus ... \oplus V_r$$
where the
$$V_i = \D^i / \D^{i-1}$$
are the quotient sheaves.  According to the remarks
above,
 $Gr(T)$ inherits the structure of a sheaf of {\sc graded nilpotent
Lie algebras}.  Here ``graded'' means that
$$[V_i, V_j ] \subset V_{i+j}.$$

If the  point $q$ is regular for $\D$ then
the  dimensions $n_i (q)$
are constant near $q$ so that the sheaves $\D^i$
correspond to  smooth subbundles of the tangent space.
In this case
$Gr(T)$ corresponds to a bundle of
Lie algebras.
(The simply connected Lie group corresponding to a particular
fiber $Gr(T_q)$ is called the {\sc nilpotentization}
of $(\D, Q)$ at $q$.)
In the Engel case
we have
$$Gr(T) = \D \oplus V_2 \oplus V_3$$
where $V_2$ and $V_3$ are real line  bundles.

Fix a nonzero element $\delta \in V_2 (q)$.
Then the Lie bracket defines a map
$\D(q) \to V_3 (q) \cong \R$,
namely $v \to [v, \delta]$.
The kernel of this map is intrinsically defined
and forms a line in the plane $\D(q)$.
Said differently,
in a neighborhood of any point where the distribution $\D$
 is Engel
there is a distinguished line field
$$L  \subset \D$$ characterized by
the fact that
$$[L, \D^2] = 0 (mod \D^2).$$
We call $L$ the {\sc Engel line field}.  It is the span of the vector
field $W$ of Engel's theorem above.

\begin{proposition}  (See \cite {Gersh}.) If an oriented
4-manifold admits an oriented Engel structure $\D$ then the manifold
is parallelizable.
\end{proposition}

{\sc Proof.} At every point $q$ of such a 4-manifold $Q$ we have the
complete flag $L(q) \subset \D(q) \subset \D^2(q) \subset T_q Q$
in the tangent bundle.  If
we show that this flag is canonically oriented
then we will be done.
(A parallelization $\{E_1,E_2, E_3, E_4 \}$
is then  be obtained by putting a Riemannian metric on $Q$
and hence on each element of the flag.  Then take
$E_1$ to be the positive unit vector spanning $L$,
$\{E_1, E_2 \}$ to be the positive orthonormal basis
for $\D$, etc.)

Write $A \in \Lambda^2 \D$ for the choice of orientation
of $\D$.    Locally $A = X \wedge Y$ for $X, Y$
nonvanishing sections of $\D$.  Then
\begin{equation}
 \delta_1 = [X,Y] \; (mod \D)
\label{eq: delta}
\end{equation}
is a well-defined
section of $T/\D$ independent of
the choices of representation X and Y.
By the assumption on the growth vector it is
non-vanishing and its span is the 
real line bundle $V_2$.  So  $\delta_1$
defines an orientation on $V_2$.  Observe that for any triple of
linear spaces, $(S, T, T/S)$ with $S\subset T$,
an orientation on any two spaces canonically determines an orientation
on the third.

 Using this observation we obtain
 an orientation on $\D^2$.
Applying the observation again we obtain an orientation
on $V_3$.  Finally, consider the
map $ad_{\delta}: \D \to V_3$ defined by
bracketing with $\delta$.
Its kernel is $L$ so that $ad_{\delta}$
induces an isomorphism
$V_3 \cong \D/L$ and so an orientation on
$\D /L$.  Applying the observation again, we finally obtain
the orientation on $L$. $\Box$

{\sc Example 1.}  If $Q$ is closed and simply connected
then it does not admit any Engel structure. This is because
its Euler class is nonzero, and hence it does not admit
a single nonvanishing vector field.
$$ $$

{\sc Example 2.}  If a 4-manifold  admits one Engel structure,
then it typically admits a continuous family of inequivalent
such structures.  For if we perturb the given structure, then we perturb
its Engel line field.  But line fields on closed manifolds typically have
continuous moduli.  For specific examples of families of
inequivalent Engel structures, see \cite
{Gersh}. The Engel situation is to be contrasted with the case of contact
structures on a three-manifold, where the moduli space  of inequivalent
structures is a discrete set.

\section{Obstruction theory and proofs.  }

Our proofs rely on the following basic
fact from topology. (See Bott and Tu,
for example.)  If $E \to  Q$
is a real oriented rank 2 vector bundle
and if $s: E \to Q$ is a section of E
{\bf which is transverse to the zero section}
$Z \subset E$ then the
section's zero locus represents the first Chern class
$c_1 (E)$ of E.  (Bott and Tu call this the Euler class.) 
Transversality implies
that this  zero locus $s^{-1} (Z) = \{q: s(q) =
0 \}$   is a canonically co-oriented smooth submanifold
of $Q$.  It represents a cohomology class $[s^{-1}(Z)]$
 via intersection theory. The value of $[s^{-1} (Z)]$
on a 2-cycle is obtained by picking a representative
for the cycle, jiggling it until it is transverse
to $s^{-1} (Z)$ and then counting intersections. 

The section $\delta_1$ of equation 1:
$$q \mapsto [X,Y](q) mod \D_q$$
is a well-defined section of $\Dp = TQ/D$
regardless of whether or not $\Sigma$ is empty.
Its zero locus is {\bf precisely}
$\Sigma_1$.  Once we have proved transversality
we will have established
\begin{lemma} 
Let $\D$ be a generic oriented rank 2 distribution   on
an oriented 4-manifold $Q$.  Then
$\Sigma_1$ is canonically co-oriented
and represents the first   Chern class
$c_1(\D^\perp)$ of the oriented rank 2 real vector bundle
$\Dp = TQ/D$. 
\end{lemma}

{\sc Proof of Lemma 1.}  It remains
to prove transversality of $\delta_1$.
Away from $C = \Sigma_1 \cap \Sigma_2$ the full growth vector
is $(2,2,4)$, meaning that $[X,[X,Y]]$
and $[Y,[X,Y]]$ (mod $\D$) span $\Dp$.
It follows that by differentiating $\delta_1$
in the $\D$-directions we obtain
all of $\Dp$.  In other words, $\delta_1$
is transverse to $Z$ at these points  and
 we only need to differentiate along
$\D$ to  achieve transversality.  

At generic points of $\Sigma_1 \cap \Sigma_2$ we will
need to differentiate in directions other than
$\D$, but transversality still holds.  
It follows from Zhitomirskii's normal form
along C ((\ref{eq:Z3}) above) that
\begin{equation}
X_1 = \partial_3
\label{eq:X1}
\end{equation}
and
\begin{equation}
X_2 = \partial_4 - x_3 ^2 \partial _1 - x_3( x_1 +f)
\partial_2
\label{eq:X2}
\end{equation}
frame $\D$ near points p of the intersection.
Here p has coordinates $(0,0,0,0)$ and 
f is a  function
satisfying $f(0) = df(0) = 0$.
We may assume that
the orientation of $\D$
is given by
$A = X_1 \wedge X_2$
since a nonzero scalar factor
will not affect transversality
considerations.  The fields
$\partial_1$ and $\partial_2$ form 
a local frame for $\Dp$. Modulo
$\D$ we have:  
\begin{equation}
\delta_1 =   -2 x_3 \partial _1 - ( x_1
+ f + x_3
\partial_3 (f) ) \partial_2  
\label{eq:d1}
\end{equation}
and
\begin{equation}
\Sigma_1 = \{ x_3 = 0 = x_1 +f \}
\label{eq:S1}
\end{equation}
near p.  One calculates at p = 0:
$ \partial_3 \delta_1 = -2 \partial_1 ,
\partial_1 \delta_1 = - \partial_2$,
mod $\D$.  Since $\partial_1, \partial_2$ span
$\D$ this proves transversality.

Finally, we need to worry about the 
codimension 4 points.  These form
a finite set, some of which may lie along
the curve $\Sigma_1 \cap \Sigma_2$.
 $\Sigma_1$ may have
singularities at these points.  But these
points do not effect cohomology. If a representative
of a  2-cycle intersects such a codimension 4 point
it can be homotoped to miss it.  And the co-orientation
is defined everywhere on $\Sigma_1$ away
from these points.  
$\Box$

To obtain the other class,  
define a section $\delta_2$ of the
vector bundle $Hom(\D, \Lambda^2 \Dp)$ by
$$\delta_2 (q) (v) = \delta_1 (q) \wedge 
([\tilde v, \tilde \delta_1 ](q) (mod \D_q)) .$$
Here $v \in \D_q$ and $\tilde v$ denotes any
extension of v to a local section of $\D$:
$\tilde v (q) =v$.  And $\tilde \delta_1$
is any local vector field with the property
that $\tilde \delta_1 = \delta_1 ( mod \D)$.
One easily checks that $\delta_2 (q)(v)$
is well-defined, independent of these
choices.  It is clear that $\D$ fails to
be Engel exactly at those points q for which
$\delta_2 (q) (v) = 0$ for all v, that is to 
say at $\delta_2 ^{-1} (0)$.  Thus
$$\Sigma := \Sigma_1 \cup \Sigma_2 = \delta_2 ^{-1} (0)$$

$\D$ and $TQ$ are oriented, so that $\Lambda^2 \Dp$
is trivial and  $Hom(\D, \Lambda^2 \Dp)
\cong \D^*$.
Thus $\delta_2$ defines a section of $\D^*$. If
this section were transverse to the zero section we would be 
done, as described at the beginning of this
section.  
However the section  $\delta_2$ {\bf cannot} be
transverse to the zero 
section at points of $C = \Si_1 \cap \Si_2$ because $\Sigma$ is
not a manifold at such points.  
In this case further analysis is needed
and this  complicates
the proof of theorem 1. The theorem
follows from 
\begin{lemma} For a generic distribution $\D$
the section $\delta_2$ is
transverse to the zero section away
from the curve $C = \Sigma_1 \cap \Sigma_2$. 
The $\delta_2$-induced co-orientations
on either $\Sigma_i$, $i =1, 2$
reverses as C is crossed.  See fig. 1. 
$C$ is homologous to zero within 
$\Sigma_1$ but might not be homologous
to zero within $\Sigma_2$.
$\Sigma$, being a Whitney stratified
set with ``cohomologically consistent co-orientations''
(see discussion below) on its principal strata, defines a
cohomology class. This class represents the first Chern
class $c_1(\D^*)$  of the dual to $\D$ bundle $\D^*$.  
\end{lemma}

\vskip 20pt
\centerline {\hbox{\hskip 0.2cm\epsfysize=6cm\epsfbox{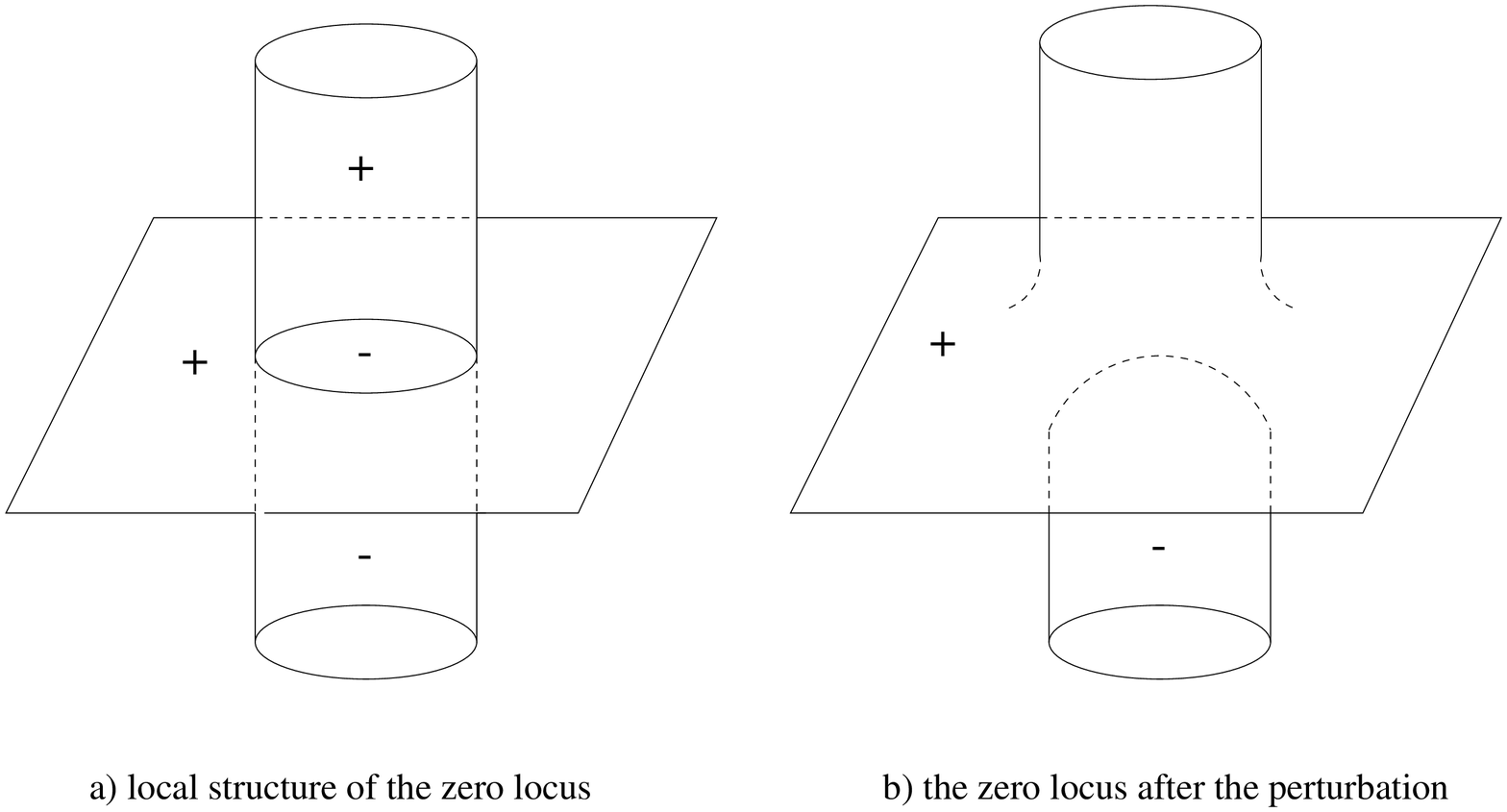}}}
$$ $$
\centerline{Fig.1. Local coorientation of $\Sigma=\Sigma_1\cup\Sigma_2$.}\vskip 5pt

We postpone the proof of the
lemma for a discussion of the
business of Whitney stratified sets
defining cohomology classes.   
it follows from Zhitomirskii's normal
forms (\ref{eq:Z1}, \ref{eq:Z2A}, \ref{eq:Z2B},
\ref{eq:Z3}, \ref{eq:Z4}) that  $\Sigma$
is a Whitney stratified subset of
$Q$ of a rather tame sort. 
Its  strata are 
the isolated codimension
4 points  $\{p_1,
\ldots, p_N \}$ on C  (their normal forms
have functional moduli and are not given by
Zhitomirskii), 
the curve $C $ minus these points  and  
the connected components of $\Sigma
\setminus C$. These latter are
 the principal strata.  They intersect
along C in a manner locally diffeomorphic to
that of two 2-planes in 4-space
intersecting along a line.  

We now   describe
how co-oriented Whitney stratified sets
may be used to define cohomology classes.  This
idea is extensively used by Vassiliev 
and we refer the reader to his book
\ {Vassiliev}, especially the introduction
and \S 8.4.  Let $Q$
be a smooth compact manifold and  $W \subset Q$
be a Whitney stratified compact subset.
Suppose that the principal strata of W
are co-oriented k-dimensional submanifolds.  
Given an $(n-k)$-cycle in Q we perturb it
slightly so that it intersects the principal strata
of W transversally and then count these intersection
points with  a plus or minus sign depending on
whether the orientation of the cycle there
agrees or disagrees with the co-orientation.  In this way
it {\em appears} that  W yields an  integer-valued function
on cycles. However 
this  number may depend on the perturbation
and may not   be
well-defined on the level of homology.  In other
words, we must somehow insure that if a cycle
Z is a  boundary then this intersection number is zero.
By cellular or simplicial approximation we
may break Z, and also the
$(n-k +1)$-cycle B which it bounds
 into a sum of
cycles each one of which is supported in an 
arbitrarily small neighborhood of Q.  First
suppose the intersection $B \cap W$
lies entirely on a principal stratum. 
Since we may assume that it is arbitrarily small,
this   intersection is diffeomorphic to
 an oriented  line segment and $Z \cap W$ consists
of two points with opposite orientations. Consequently
its intersection numbers add to zero as desired.
If $B$ intersects a stratum R  of dimension $k-1$
we must impose a {\sc cohomological consistency
condition}, to be spelled out momentarily,  on
W near R.  If B intersects
a stratum of dimension $k-2$ or less then we may slightly
perturb it while keeping its boundary
fixed, and in this way  reduce to the case
of intersection with a strata of dimension $k$ or $k-1$.

{\sc Cohomological Consistency Condition:} The link L of
any stratum of dimension $k-1$ is cohomologically trivial
within the sphere of dimension $n-k$. {\bf See figure 2.}

This link is a finite collection of co-oriented
points
on the $(n-k)$-sphere.  The condition is then that
the sum of co-orientations is zero.
In particular the number of points in
the link must be even.  

{\sc Our case:}  {\b 
See the right hand picture of figure 2.}
In our case $k =2$ so
$k-1 = 1$.  The 1-dimensional stratum 
consists of the smooth points of  C,
that is all the points of C except the
points of codimension 4 not covered in
Zhitomirskii's paper.
The link of C consists of 4 points,
representing the 4 `pieces' of 
$\Sigma \setminus C$ near C.  The cohomological
condition is that two of these
points have plus signs and two have 
minus signs. 

\vskip 20pt
\centerline {\hbox{\hskip 0.2cm\epsfysize=3.5cm\epsfbox{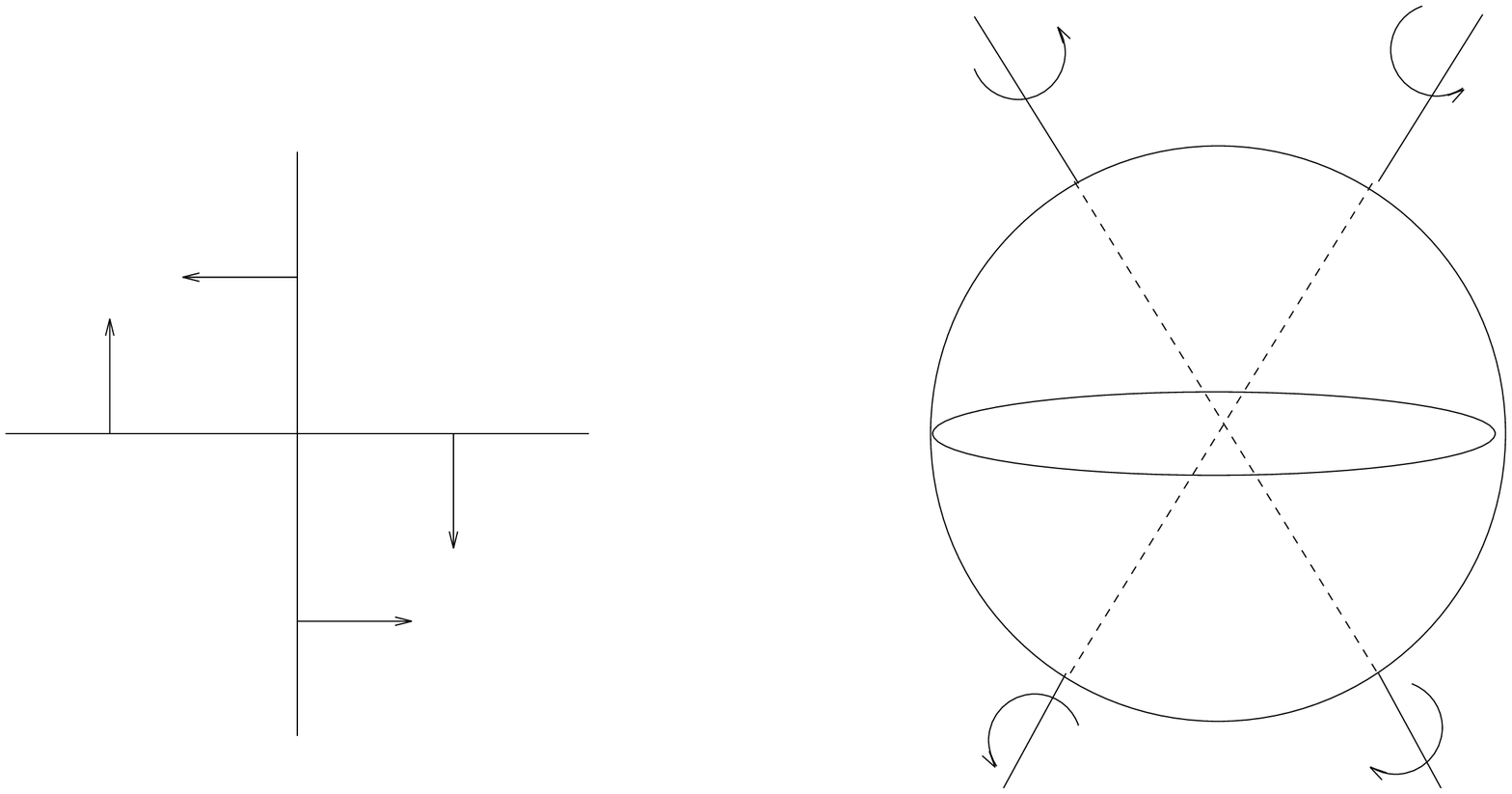}}}
$$ $$
\centerline{Fig.2. Local coorientations for the case $n-k=1$ and $n-k=2$.}\vskip 5pt 

We recall that   the notion of the link L 
of a singular stratum R. R is a piece of
a submanifold of dimension r, where  
 $r \le
k-1$. Intersect R
transversally at a point p with  a small piece of
an $(n-r)$-manifold V. Intersect the  result with an
$(n-r -1)$-sphere
$S \subset V$  surrounding p.
We call L, or the pair $(L, S)$ the {\sc link} of
R at p.  By the definition of
a stratified set, $V \cap W$ is
diffeomorphic to the cone $C(L)$  over L, and
a  p admits a neighborhood N and
a diffeomorphism which takes the pair
$(W \cap N, N)$ to $(U \times C(L), U \times V)$ for
some open set U in R.   

\begin{lemma} The cohomological condition
insures that W is a well-defined cohomology
class.
\end{lemma}

{\sc Proof.}  It remains to show 
that if the small $(n-k+1)$-cycle B intersects
a stratum of dimension $k-1$ then the
intersection number of its boundary Z
with W is zero. We may take the cycle B 
to be ball intersecting R transversally.
Its boundary can be homotoped
so as to form the sphere S used to 
define the link L of R.  Then the intersection
$Z \cap W$ is the sum of the points
L with appropriate intersection numbers,
which we have assumed to be zero.    
$\Box$

{\sc proof of lemma 2.}
 If $\delta_2$
is transverse to the zero section
away from C as claimed, and if the
induced co-orientations on $\Si \setminus C$
 reverse 
as C is crossed while travelling along a fixed
$\Sigma_i$, then the cohomological 
consistency condition follows directly.   For
the link consists of four points and this
reversal implies they come in pairs which cancel,
one pair for each $\Sigma_i$.  See figure 2.

The claim regarding the Chern class
also follows directly.  To see this
imagine a cycle K and its intersection
number with $\Sigma$ as we have just defined it.
The intersection points with K
lie outside some small neighborhood
U of C.  By Sard's theorem,
we may perturb the section $\delta_2$
so that the resulting section s is transverse
to the zero section Z of $\D^*$.
Moreover, this perturbation may be concentrated
within U
so that s agrees with $\delta_2$ outside 
of U.  Then the intersection
number of K with $\Sigma$ equals that of K with
$s^{-1}(0)$.  The later represents the first Chern
class of $\D^*$ since s is a transverse section.   
{\bf See figure 1 again.}

It remains to check the transversality claims.
We will only check them near C.  For
points of $\Si$ away from C or L
transverality follows from the normal
forms  (\ref{eq:Z1}),  (\ref{eq:Z2A}),
(\ref{eq:Z2B})  
given above. The calculation follows
the lines we follow below, but are 
significantly simpler.  For points along
L  the calculation is similar  to thaat below
and is also omitted. We would use the
form (\ref{eq:Z4}) .

For points along C we use (\ref{eq:Z1})  in
the dual form of  
 equations (\ref{eq:X1}) - (\ref{eq:S1}) 
which follow lemma 1.
We have
$$[X_1 , \delta_1] = -2 \partial_1 - (2 \partial_3 f
+ x_3 \partial_3 ^2 f) \partial_2
 $$ and
$$[X_2, \delta_2] = \{X_2 [f + x_3 \partial_3 f] -
\delta_1 [x_3(x_1 + f)] 
\}
\partial_2
= \{ \partial_4 f + O(2) \} \partial_2 
$$ 

If $x \in \Sigma_1$ then $\delta_1$ must be
proportional to $[X_1, \delta_1]$.
Multiplying our expression for $[X_1, \delta_1]$
by $2 x_3$ and refering
back to equation (\ref{eq:d1}) 
for $\delta_1$ we see that this means that:
\begin{equation}
x_1 +f = x_3 h
\label{eq:S21}
\end{equation}
where h is a certain function vanishing to 1st order.
Also we must have $[X_2, \delta_1] = 0$
when $x_3 \ne 0$ since in this case
$\delta_1$ has a $\partial_1$-component
whereas $[X_2, \delta_1]$ has only a   $\partial_2$-component.  
Now we may assume, by genericity,
that
$d(\partial_4 f)   \ne 0$ so that the
leading term of this $\partial_2$-component
is $\partial_4 f$.  We may then
write $[X_2, \delta_2] = {\bar x_4} \partial_2$
where $\bar x_4 = \partial_4 f + O(2)$ is
a new coordinate.  Also set
$$  \tilde x_1  = x_1 + f.$$
 
Generically we have $d \tilde x_1  \wedge dx_2
\wedge dx_3
\wedge d \bar x_4 \ne 0$
so that $( \tilde x_1, x_2 , x_3 , \bar x_4)$
form coordinates near p. Observe that
$\Sigma_1$ can equally well be expressed
as the locus of points with
$x_3 = 0$ and $\tilde x_1 - x_3 h = 0$.
Set
$$
\bar x_1 = \tilde x_1 - x_3 h$$ 
Thus locally
$$ \Sigma_1 = \{ x_3 = 0, \bar x_1 = 0\} ,$$
$$ \Sigma_2 = \{ \bar x_1 = 0 , \bar x_4 = 0\} ,$$
and
$$ C = \{ 
 x_3 = 0, \bar x_1 = 0 , \bar x_4 = 0\} .$$
Note that this expresses the intersection C
as a line obtained by intersecting two
2-planes in 4-space as claimed.  

The section $\delta_2$ is a fiber-linear map which
is given on the basis $\{X_1, X_2\}$ by
$X_i \mapsto \delta_1 \wedge [X_i, \delta_1]$.
Using the above formulae we
compute that 
we can represent $\delta_2$ as the
 2-vector:
\begin{equation}
\delta_2 = ( - 2 \bar x_1 + x_3 g,  -2 x_3 \bar x_4)
\label{eq:d2}
\end{equation}
where g is function vanishing to 1st order at p.
We have also dropped off the basis vector coefficient
$\partial_1 \wedge 
\partial_2 
$ from $\delta_2$.  (Here $(1,0), (0,1)$ represent the dual basis
to $X_1, X_2$.)
Then near p we have:
$${\partial \over {\partial \bar x_1}} \delta_2
= ( -2, 0 )  ,$$
$${\partial \over {\partial  x_3}} \delta_2
= ( g, -2 x_4 ) .$$
If $X_1 \wedge X_2$ represents the orientation
of $\D$ then $ \{ (1,0)    , 
(0,1)  \}$
represent the corresponding oriented frame for
$Hom(\D, \Lambda^2 \D) = \D^*$.
Observe  that the normal bundle to
$\Sigma_1$ near p is framed by 
${\partial \over {\partial \bar x_1}} $
and 
${\partial \over {\partial  x_3}}$.  
It follows from our expressions for the
derivatives of $\delta_2$ 
that the $\delta_2$-induced co-orientation of $\Sigma_1
\setminus C$ is given by the frame
$\{ {\partial \over {\partial \bar x_1}},
 \bar x_4 {\partial \over {\partial  x_3}} \}$.
 The induced co-orientation reverses
as we cross C travelling on $\Sigma_1$,
since $C$ is defined on $\Sigma_1$
by $\bar x_4 = 0 $. 
Also, the section $\delta_2$ is transverse
to the zero section for points of 
$\Sigma_1$ away from C.
 
To complete the proof we compute:
$${\partial \over {\partial  \bar x_4}} \delta_2
= ( 0 , -2 x_3 ) .$$
The vector fields 
${\partial \over {\partial  \bar x_1}} $ and
${\partial \over {\partial  \bar x_4}} $
frame the normal bundle to $\Sigma_2$ near p.
For the same reasons as above we see 
that the $\delta_2$-induced co-orientation
of $\Sigma_2$ is given by the frame 
$\{ {\partial \over {\partial \bar x_1}},
 x_3 {\partial \over {\partial  x_4}} \}$.
On $\Sigma_2$ the curve C is given by
$x_3 = 0$.  Again the induced co-orientation
reverses as we cross $C$ and the section
is transverse away from C.
$\Box$

{\sc Remark.} The co-orientability of $\Si_2$ depends on the orientability
of the $3$-distribution $D^3\vert_{\Si_2}$.
One can check 
that the subvariety in the $2$-jets of germs of $2$-distributions
corresponding to the
growth vector $(2,3,3)$ is not coorientable.

{\sc Proof of Theorem 2.} The construction of the sections $\delta_1$ and
$\delta_2$
can be adjusted to the nonoriented case to define the sections of bundles
$Hom (\Lambda^2 \D,\D^\perp)$ and $Hom (\D,\Lambda^2(\D^\perp))$
such that $\Si_1$ and $\Si=\Si_1\cup\Si_2$ are the zero loci of
these sections. Therefore, these cycles are dual to the
Stiefel-Whitney classes of the corresponding bundles which can
be calculated by using the standard methods of the theory of
characteristic classes. In particular, we have,
$$[\Si_1]=w_2(Hom
(\Lambda^2 \D,\D^\perp))=w_1^2(\D)+w_2(\D)+w_2(Q);$$
$$[\Si_1\cup\Si_2]=w_2(Hom
(\D,\Lambda^2(\D^\perp))=w_2(\D)+w_1^2(Q)+w_1(\D)w_1(Q);$$
$$[\Si_2]=[\Si_1]+[\Si_1\cup\Si_2]=
w_1^2(\D)+w_1(Q)w_1(\D)+w_1^2(Q)+w_2(Q).$$

If the restriction of $\D$ to $\Si_1$ is orientable then the arguments at
the beginning of this
section  show that $\Si_1\cap\Si_2$ is a boundary. Therefore
$\Si_1\cap\Si_2$ represents the  homology class on $\Si_1$ dual to
$w_1(\D \vert_{\Si_1}).$ (Using a bundle 
homomorphism which  is essentially
given by $\delta_2$ one can explicitly find a bundle over $\Si_2$
isomorphic  to
$\Lambda^2 \D \vert_{\Si_1}$ and a section of it for which $\Si_1\cap
\Si_2$ is the zero locus).
Hence, we can apply the Gysin formula to compute the cohomology class dual  to
$\Si_1\cap\Si_2$ as follows
$$[\Si_1\cap\Si_2]=i_*([\Si_1\cap \Si_2\vert_{\Si_1}])=
i_*(w_1(\D\vert_{\Si_1}))=w_1(\D)i_*(1)=w_1(\D)[\Si_1].$$
     Here $[\Si_1\cap \Si_2\vert_{\Si_1}]$ denotes the class dual to the
cycle $\Si_1\cap
\Si_2$ in the cohomology group of $\Si_1$, $i:\Si_1\to Q$ is the natural
inclusion and $i_*$
is  the corresponding Gysin homomorphism in cohomology. $\Box$

\section{Open problems.}

{\sc Problem 1.} Does every closed parallelizable 4-manifold admit an
Engel structure?  (The proof suggested in \cite {Gersh} is incomplete.)

This problem, the converse to the parallelizability proposition,
is the basic open question in the area.  The  vanishing
of the obstructions $c_1(\D^\perp)$ and $c_1(\D^*)$ of  our main  theorem
implies the parallelizability of $Q$. On the other hand, a
parallelizable  $Q$ typically admits countably many homotopically
distinct 2-distributions $\D$ satisfying these  vanishing conditions.
The  next problem is a stronger version of problem 1.

{\sc Problem 2.} Suppose that the rank two distribution $\D$ satisfies
$c_1(\D^\perp)= c_1(\D^*)=0$. Can we homotope $\D$ to an Engel
distribution?

\end{document}